\documentclass[preprint,review,12pt]{elsarticle}

\usepackage{amssymb}
\begin{document}

\begin{frontmatter}



\title{Extended maximum
concurrent flow problem with saturated capacity}


\author{Congdian Cheng}
%\footnotetext
%{
%\baselineskip 10pt Received -- --, 2005; accepted
%-- --, 2005\\  DOI:
%10.1007/s11425-005-0036-y\\
%Corresponding author: Congdian Cheng.\\ E-mail: zhiyang918{\rm @}163.com.\\
%Post address: College of Mathematics and Systems Science, Shenyang
%Normal University,
%  Shenyang,  110034,  People's Republic of China.}
%\footnotetext

\address{ College of Mathematics and Systems Science, Shenyang Normal University,
  Shenyang 110034, PR China\\ E-mail: zhiyang918{\rm @}163.com.}

\begin{abstract}
The present work studies a kind of  Maximum Concurrent Flow Problem,
called as Extended Maximum Concurrent Flow Problem with Saturated
Capacity. Our major contributions are as follows: (A) Propose the
definition of Extensive Maximum Concurrent Flow Problem with
Saturated Capacity  and prove its solutions exist. (B) Design a
approximation algorithm to solve the problem. (C) Propose and prove
the complexity and  the approximation
 measures of the  algorithm we design.
\end{abstract}


\end{frontmatter}

% \linenumbers

\section{Introduction}
\label{}Maximum Concurrent Flow Problem (MCFP) is a main kind of
multicommodity flow problems, which was introduced by Matula in
1985, see Shahrokhi and  Matula [9],  and have been severely studied
for more than two decades, see, e.g., [2,4,5,6,9,10]. Motivated by
the published literatures of this research area, we discuss a kind
of specific Maximum Concurrent Flow Problem in the present work,
termed as  Extended Maximum Concurrent Flow Problem with Saturated
Capacity (EMCFPSC).

The rest of this article is organized as follows. Some preliminaries
are  presented in Section 2. Section 3 formulates the  problem
EMCFPSC and prove the existence of its solutions. Section 4 is
specially devoted to designing an algorithm to approximately solve
the problem. Section 5 investigates the complexity and precision of
the algorithm. Finally, the paper is  concluded with Section 6.

\section{Preliminaries}
This section provides some preliminaries for our sequel research.

A graph ${G=(V,E)}$ is called as a hybrid graph
%\footnote{For the
%basic knowledge of digraph, please see [8] or other references.}
 if
its edge can be either directed edge or undirected edge, of which
the directed graph and undirected graph are specific cases. Here
 $V$ and $E$ represent all the nodes (or vertices) and edges (or
arcs) of $G$ respectively.

Given graph ${G=(V,E)}$  with weight (or capacity) $c$ and
$H=\{[s_{i}, t_{i}]:s_{i},t_{i}\in V; i=1,2,\cdots,k\}$,  where
$c:E\longrightarrow R_{+}$ and $R_{+}=[0,\infty)$, and $s_{i}$ and
$t_{i}$ express the source and the terminal of commodity $i$
respectively,  we call the triad $(G,c,H)$ as a multicommodity
network. Let ${[s_{i},t_{i}]\in H}$;
$s_{i},v_{1},v_{2},\cdots,v_{l},t_{i}$ be some nodes of $G$ and be
different from each other except for $s_{i}=t_{i}$; $e_{j}$ be the
edge of $G$ with endpoint $v_{j}$ resp. $v_{j+1}$ for $0\leq j\leq
l$, where $v_{0}=s_{i}$ and $v_{l+1}=t_{i}$. (When $e_{j}$ is
directed,
 the edges $v_{j}$ and $v_{j+1}$
must be the original endpoint
 and  the terminal endpoint, respectively.) Then we call
 $P=[(s_{i},e_{0},v_{1}),(v_{1},e_{1},v_{2}),
\cdots,(v_{l},e_{l},t_{i})]$ as an ${s_{i}-t_{i}-}$ path of
$(G,c,H)$. (To simplify in notation afterwards, we denote
$(v_{j},e_{j},v_{j+1})$ as $(v_{j},v_{j+1})$ when $e_{j}$ needn't be
indicated; and denote $(v_{j},e_{j},v_{j+1})$ as $e_{j}$ when
$v_{j}$ and $v_{j+1}$ needn't  be indicated.) Let
${\mathcal{P}}_{i}$ be a set of ${s_{i}-t_{i}-}$ paths such that
$c(e)>0$ for all $e\in E(P)$ and $P\in{\mathcal{P}}_{i}$,
$i=1,2,\cdots,k$, where $E(P)$ denotes  the set of all the edges of
$P$. We call ${\mathcal{P}}=\bigcup\limits_{i=1}^{k}
{\mathcal{P}}_{i}$ as a (positive) path system  on $(G,c,H)$, which
is denoted by ${\mathcal{P}}|_{(G,c,H)}$, or ${\mathcal{P}}$ for
conciseness.

\noindent{\bf{Definition 1.}} Let ${\mathcal{P}}$ be a path system.
If mapping ${y:{\mathcal{P}}\longrightarrow R_{+}}$ satisfies:
$$\sum\limits_{e\in E(P),P\in {\mathcal{P}}} y(P)\leq
c(e)\hspace*{4mm}\forall e\in E({\mathcal{P}}),$$ where
$E({\mathcal{P}})=\{e\in E(P):P\in {\mathcal{P}}\}$, i.e. the set of
all the edges of ${\mathcal{P}}$. Then we call ${y}$ as a
multicommodity flow on ${\mathcal{P}}$ defined by the function on
 the path system, and as a flow for simplicity; and call ${V(y)=\sum\limits_{P\in
{\mathcal{P}}}y(P)}$ as flow value of $y$, and
${V_{i}(y)=\sum\limits_{P\in {\mathcal{P}}_{i}}y(P)}$ as branch flow
value of commodity $i$ with $y$. The set of all the flows on
${\mathcal{P}}$ is denoted by $F[{\mathcal{P}}]$, and by $F$ in
brief.

\noindent{\bf Remark 1.} {As is known to us, there are two kinds of
definitions about the flow, which are the defined by the function on
the edge set and by the function on the path system. For the purpose
of being simple and clear, we only make the study as far as the flow
defined by the function on the path system in the present work.}

Generally, the problem to find a multicommodity flow to satisfy
certain conditions is called the multicommodity flow problem (MFP).
A few of usual multicommodity flow problems are as follows.  The
problem to find a flow $y$  so that $V(y)$ is maximized, i.e.
$V(y)=\textrm{OPT}[{\mathcal{P}}]=\hat{V}(=max\{V(y):y\in F[P] \})$,
is called the maximum multicommodity flow problem on ${\mathcal{P}}$
  (MMFP). A flow $y$ satisfying MMFP is called a
solution of MMFP.  The problem to find a flow $y$  so that
$V_{i}(y)\leq b_{i}$ for given $b_{i}$, $i=1,2,\cdots,k$, and $V(y)$
is maximized, i.e. $V(y)=\overline{V}(=max\{V(y):V_i(y)\leq
b_i,i=1,2,\cdots,k;y\in F\})$, is called the
 maximum flow problem with bounds(MMFP-B).
The problem to find a flow $y$  so that
 $V(y)$ is maximized
under the condition $V_{i}(y)= \lambda b_{i}$ for given $b_{i}$,
$i=1,2,\cdots,k$, and for $\lambda\in[0,1]$ is called the Maximum
Concurrent Flow Problem (MCFP) (see, e.g., Shahrokhi and  Matula
[10], or Garg and K\"{o}nemann [6]).

\section{Model formulation}
Now we begin to introduce the problem we are attacking in the
present work.

\noindent {\bf {Definition 2.}} Suppose $\mathcal{P}$ is a path
system on the multicommodity network $(G,c,H)$, and
$\mathbf{b}=(b_{1},b_{2},\cdots,b_{k})$.  We call the problem to
find a flow $y'$ so that $\min\limits_ { 1\leq i\leq k
}[\frac{1}{b_{i}}V_{i}(y')] =\max\{\min\limits_ { 1\leq i\leq k
}[\frac{1}{b_{i}}V_{i}(y)]:y\in F, V_{i}(y)\leq b_{i}\}$ as Extended
Maximum Concurrent Flow Problem (EMCFP). We call the problem to find
a flow $y'$ so that $\min\limits_ { 1\leq i\leq k
}[\frac{1}{b_{i}}V_{i}(y')] =\max\{\min\limits_ { 1\leq i\leq k
}[\frac{1}{b_{i}}V_{i}(y)]: y\in F, V_{i}(y)\leq b_{i}\}$, and
$V(y')=\max\{V(y): y$ is the solution of EMCFP $\}$ as Extended
Maximum Concurrent Flow Problem with Saturated Capacity(EMCFPSC).

\noindent{\bf Theorem 1.} The solution of the problem EMCFPSC
exists.

\noindent{\bf\it Proof.} It is obvious that the solution of MCFP is
also the solution of EMCFP. Hence,  $B=\{y\in F: y$ is the solution
of $FMCFP\}\neq\emptyset$ from the well known fact the solution of
the problem MCFP exists. When the member of $B$ is finite, the
conclusion of  Theorem 1 is trivial. Otherwise, $\sup\{V(y): y\in
B\}$ exists and is finite. Let it be $b$. Then there is at least a
sequence  $\{y_j\}$ in $B$ such that
$\lim\limits_{j\rightarrow\infty}[\sum\limits_{i=1}^k V(y)=b$, and
$\lim\limits_{j\rightarrow\infty}y_j(P)$ exists for all $P\in
{\mathcal{P}}$. Now let
 $y(P)=\lim\limits_{j\rightarrow\infty}y_j(P)$ for all $P\in
{\mathcal{P}}$. Then $y=\in B$, and $V(y)=b$. That is, $y$ is the
solution of EMCFPSC.\hspace{1cm}$\Box$

\noindent{\bf Example} During the winter between the year 2009 and
2010, in major regions of the world, the amount of energy consumed
for heating is far more than the ordinary winter due to the abnormal
lower air temperatures. With this circumstances, the electrical
networks emerge to congestion for a lot of districts in the world.
Assume the purpose of supplying power is firstly to maximize the
minimum of the met rates of  cities. Then the problem  to optimize
the power supplement can be largely  regarded as the mould EMCFPSC.
One can believe the  mean of applications of the mould EMCFPSC from
the instance.

\section{Algorithm}
In this section, we are devoted to designing an approximation
algorithm for the problem EMCFPSC.

\noindent{\bf Algorithm} (Approximation algorithm for EMCFPSC):

\noindent{\bf Input} Multicommodity network $(G, c, H)$, path system
$\mathcal{P}$, vector $\mathbf{b}$ and error parameter $0<\eta<1$.

\noindent{\bf Output} Flow $y^\ast$ on $\mathcal{P}$, which has the
characters specified by the Theorem in Section 5.

\begin{enumerate}\item Put $\epsilon=\min\{\frac{\eta}{\sum\limits_{i=1}^{k}b_i},
\frac{1}{2}\}; b_i^0= b_i, i=1,2,\cdots,k; l=0, h=0$.
\item Set $l=l+1$. If $l\eta>1$, implement the next Step.
Otherwise, put $l=l^*$ and go to the Step 4.
\item Set $b_i=l\eta b_i^0, \mathbf{b}=(b_1,b_2,\cdots,b_n)$.
For $\mathcal{P}, \mathbf{b}$ and $\epsilon$, obtain an
approximation solution $y_{1}$ of MMFP-B with by Algorithm 2 of
Cheng [3]. If $V(y_{1})< \frac{1}{1+\epsilon}
 (\sum\limits_{i=1}^{k}b_i)$, put $l=l^*$ and implement the next step. Otherwise,
 return to (2).
\item Construct the auxiliary network and demanding vector $\mathbf{b}$ as follows.
Set $V'=V\bigcup\{t'_0,t'_i:i=1,2,
\cdots,k\}(t'_0,t'_i\bar{\in}V),E'=E\bigcup\{(t_i,t'_0),
(t_i,t'_i):i=1,2,\cdots,k\},
H=\{[s_i,t'_i],[s_i,t'_0]:i=1,2,\cdots,k\}, G'=(V', E');
{\mathcal{P'}}_i=\{P+(t_i,t'_i):P\in {\mathcal{P}}_i\},
{{\mathcal{P}}'}_0^i=\{P+(t_i,t'_0):P\in
{\mathcal{P}}_i\},i=1,2,\cdots,k;
{\mathcal{P}}'=(\bigcup\limits_{i=1}^k{\mathcal{P}}'_i)
\bigcup(\bigcup\limits_{i=1}^k{{\mathcal{P}}'}_0^i);
b'_i=(l^*-1)\eta b_i^0, b'_0=b,
\mathbf{b'}=(b'_1,b'_2,\cdots,b'_k,b'_0)$, and
$$c'(e)=\left\{\begin{array}{ll}
c(e)\hspace*{9 mm}e\in E(G)\\ b'_i\hspace*{12mm}e=(t_i,t'_i)\\
b_i-b'_i \hspace*{10mm}e=(t_i,t^0_i)\end{array}\right..$$
\item Put $h:=h+1, b=h\eta$. If $b> \sum\limits_{i=1}^{k}b_i^0$,
put $h^*=h$ and implement the next step. Otherwise
 for ${\mathcal{P}}', b'$ and $\epsilon$, obtain an
 approximation solution $y_2$
 by Algorithm 2 of Cheng [3]. If $V(y_2)<
\frac{1}{1+\epsilon}
 [(\sum\limits_{i=1}^{k}b'_i)+b]$, put $h^*=h$ and implement the next step.
 Otherwise, put $y'=y_2$, and then return  to 5.
\item Put $y^\ast(P)=y'(P+(t_i,t'_i))+
y'(P+(t_i,t_0)),\forall P\in {\mathcal{P}_i}, i=1, 2, \cdots, k$.
Stop.
\end{enumerate}

\section{Algorithm analysis}
This section is specially devoted to discussing the correctness,
approximate
    precision and complexity of the algorithms we nave presented above.

\noindent{\bf Theorem} The complexity of the Algorithm  is
$O(\frac{1}{\eta^3}km^2\log n)$, where $k=|H|,n=|V(G)|,m=|E(G)|$.
Let $y$ be a solution of EMCFPSC and $y^\ast$ be the output of the
Algorithm. Then $y^\ast$ is a flow on ${\mathcal{P}}$,
$V_i(y^\ast)\leq b_i, i=1,2,\cdots,k,
V(y)\leq[l^*(\sum\limits_{i=1}^{k}b_i)+h^*]\eta$, and
$$[(l^*-1)(\sum\limits_{i=1}^{k}
b_i) +(h^*-1)]\eta-2\eta\leq
V(y^\ast)\leq[(l^*-1)(\sum\limits_{i=1}^{k}b_i)+h^*]\eta.$$
$$(l^*-1)\eta
-\frac{2\eta}{\min
b_i}\leq\min\frac{V_i(y^\ast)}{b_i}\leq\min\frac{V_i(y)}{b_i}\leq
l^*\eta.$$

\noindent{\bf\it Proof.} It is obvious that the complexity of the
Algorithm depends on the complexity of  the  subroutine Algorithm 2
of [7] and the times of the iteration that we operate on  the
subroutine. By simply analyzing the Algorithm, we can know that the
times  is no more than $\max b_i(\frac{k}{\eta}+\frac{1}{\eta})$,
which can be denoted as $O(\frac{1}{\eta})$. On the other hand, the
complexity of the subroutine  is $O(\frac{1}{\epsilon^2}km^2\log
n)=O(\frac{1}{\eta^2}km^2\log n)$, see [3]. Hence the complexity of
Algorithm  is $O((\frac{1}{\eta^3}km^2\log n)$.

By further analyzing the Algorithm, we can easily know that $y^\ast$
is a flow on $\mathcal{P}$, $V_i(y^\ast)\leq b_i, i=1,2,\cdots,k$,
$V(y)\leq[l^*(\sum\limits_{i=1}^{k}b_i)+h^*]\eta,
V(y^\ast)\leq[(l^*-1)(\sum\limits_{i=1}^{k}b_i)+h^*]\eta$, and
$\min\frac{V_i(y^\ast)}{b_i}\leq\min\frac{V_i(y)}{b_i}\leq l^*\eta$.
Hence we only specify $[(l^*-1)(\sum\limits_{i=1}^{k} b_i)
+(h^*-1)]\eta-2\eta\leq V(y^*)$ and
$\min\frac{V_i(y^*)}{b_i}\geq(l^*-1)\eta -\frac{2\eta}{\min b_i}$.

In terms of  the Algorithm, we have $\epsilon\leq
\frac{\eta}{\sum\limits_{i=1}^{k} b_i}, (l^*-1)\eta\leq 1,
(h^*-1)\eta\leq \sum\limits_{i=1}^{k} b_i$, and $V(y^*)\geq
\frac{1}{1+\epsilon}[(\sum\limits_{i=1}^{k}(l^*-1)\eta b_i)
+(h^*-1)\eta]$. Therefore
$$\begin{array}{lcl} V(y^*)&\geq& \frac{1}{1+\epsilon}
[(\sum\limits_{i=1}^{k}(l^*-1)\eta b_i)
+(h^*-1)\eta]\\
&=&[(l^*-1)(\sum\limits_{i=1}^{k} b_i)
+(h^*-1)]\eta-\frac{\epsilon}{1+\epsilon}
[(l^*-1)\eta(\sum\limits_{i=1}^{k} b_i)\\& & +(h^*-1)\eta]\\
&\geq& [(l^*-1)(\sum\limits_{i=1}^{k} b_i)
+(h^*-1)]\eta-\epsilon[(\sum\limits_{i=1}^{k} b_i)
+(\sum\limits_{i=1}^{k} b_i)]\\&\geq& [(l^*-1)(\sum\limits_{i=1}^{k}
b_i) +(h^*-1)]\eta-2\eta.
\end{array}\eqno(1)$$

Let $V'_i=\sum\limits_{P\in{\mathcal{P}}'_i} y'_i(P), i=1, 2,
\cdots, k$ and $V'^0_i=\sum\limits_{P\in{\mathcal{P}}'^0_i}
y'^0_i(P)$. Then $V_i(y^*)=V'_i+V'^0_i$ since $y_i=y'_i+y'^0_i$. As
$b'_i=(l^*-1)\eta b_i$, we can have $V'_i\leq (l^*-1)\eta b_i$.
Suppose $V_j(y^*)<(l^*-1)\eta b_j-2\eta$ for some $j$. Then
$$\begin{array}{lcl} V(y^*)&=&\sum\limits V_i=[\sum\limits_{i\neq j}
(V'_i+V'^0_i)]+V_j(y^*)\\
&<&\sum\limits_{i\neq j} V'_i+\sum\limits_{i\neq
j}V'^0_i+(l^*-1)\eta b_j-2\eta.\end{array}\eqno(2)$$ On the other
hand, $\sum\limits_{i\neq j}V'^0_i\leq\sum V'^0_i\leq(h^*-1)\eta$.
Hence, (2) implies
$$\begin{array}{lcl}V(y^*) &<&[(\sum\limits_{i\neq
j}
(l^*-1)\eta b_j)+(h^*-1)\eta]+(l^*-1)\eta b_j-2\eta\\
&=&[(\sum\limits (l^*-1)\eta b_j)+(h^*-1)\eta]-2\eta.
\end{array}\eqno(3)$$
It  is clearly that (3) is in contradiction with (1). So,
$V_i(y^*)\geq (l^*-1)\eta b_i-2\eta$ for any $i$. This implies
$$\min\frac{V_i(y^*)}{b_i}\geq(l^*-1)\eta -\frac{2\eta}{\min b_i}.$$
This comples the proof.\hspace{1cm}$\Box$

\noindent{\bf Remark 2.} Recently,  B\"{u}sing and  Stiller [1]
considered a kind of network flow problem arising in line planning.
It may be an interesting topic to explore the application of the
present work in line planning. Moreover, Soleimani-damaneh [11]
investigated the determination of maximal flow in a fuzzy dynamic
network with multiple sinks.  Mehri [8] studied the inverse maximum
dynamic flow problem. We believe to extend the present work in  a
fuzzy dynamic network and  to address the  inverse problem EMCFPSC
may also be  interesting topics.
\section{Concluding remarks}
In this paper, the problem Extended Maximum Concurrent Flow Problem
with Saturated Capacity is introduced. The existence of its
solutions is proved.  A approximation algorithm to solve the problem
is designed, and the effectiveness of the algorithm, including the
complexity and approximation
 measures,  is discussed. The approach we design the algorithm is
 the specific contribution of the present work, whose main
 characters
 are   to construct auxiliary network  and to implement iterating search
 through subroutine a known algorithm. To solve other network problems with this approach
 is an interesting
 topic for further researches in the future.
\vskip .5cm
   \noindent{{\bf Acknowledgements}
The authors cordially thank the anonymous referees
    for their valuable comments which lead to the improvement of
    this paper.

This research is Supported by the Eduction Department Foundation of
Liaoning Province (L2010514).
   }



%% References


%%\bibliographystyle{model1-num-names}
%%\bibliography{<your-bib-database>}



\end{document}